 \documentclass[12pt]{article}

 \usepackage{mathrsfs}
 \usepackage{amssymb}
 \usepackage{color}

 \usepackage{amsfonts,amsmath,mathrsfs,amssymb}
 \usepackage{times,helvet,courier,type1cm}
 \usepackage{harvard}

 \allowdisplaybreaks

 \newtheorem{thm0}{Theorem}[section]
 \newtheorem{exa0}{Theorem}[section]

 \newtheorem{def1}[thm0]{Definition}
 \newtheorem{lem1}[thm0]{Lemma}
 \newtheorem{thm1}[thm0]{Theorem}
 \newtheorem{cor1}[thm0]{Corollary}
 \newtheorem{pro1}[thm0]{Proposition}
 \newtheorem{con1}[thm0]{Condition}
 \newtheorem{rem1}[thm0]{Remark}
 \newtheorem{exa1}[exa0]{\it{Example}}

 \def\bglemma{\begin{lem1}}\def\edlemma{\end{lem1}}
 \def\bgtheorem{\begin{thm1}}\def\edtheorem{\end{thm1}}
 
 \def\bgproposition{\begin{pro1}}\def\edproposition{\end{pro1}}

 \def\benumerate{\begin{enumerate}}\def\eenumerate{\end{enumerate}}
 \def\bitemize{\begin{itemize}}\def\eitemize{\end{itemize}}

 \def\beqlb{\begin{eqnarray}}\def\eeqlb{\end{eqnarray}}
 \def\beqnn{\begin{eqnarray*}}\def\eeqnn{\end{eqnarray*}}

 \def\eqref#1{{\rm(\ref{#1})}}

 \def\<{\langle}\def\>{\rangle}

 \def\mbb{\mathbb}\def\mbf{\mathbf}
 \def\mrm{\mathrm}

 \def\proof{\noindent{\textit{Proof.~~}}}
 \def\qed{\hfill$\square$\smallskip}
 \def\beqlb{\begin{eqnarray}}\def\eeqlb{\end{eqnarray}}
 \def\beqnn{\begin{eqnarray*}}\def\eeqnn{\end{eqnarray*}}

 \def\qqquad{\qquad\qquad}

 \def\beqlb{\begin{eqnarray}}\def\eeqlb{\end{eqnarray}}
 \def\beqnn{\begin{eqnarray*}}\def\eeqnn{\end{eqnarray*}}

 \def\mbb{\mathbb}

 \def\qqquad{\qquad\qquad}
 \def\qed{\hfill$\Box$\medskip}

 \def\<{\langle}\def\>{\rangle}

 \def\mbb{\mathbb}\def\mrm{\mathrm}\def\mbf{\mathbf}

\begin{document}

\bigskip
\centerline{\bf\LARGE Ergodic and strong Feller properties}

\medskip

\centerline{\bf\LARGE of affine processes\footnote{ Supported by the National Natural Science Foundation of China (No. 11531001).}}

\bigskip

\centerline{Shukai Chen and Zenghu Li}

\medskip

\centerline{\it Department of Mathematics, Beijing Normal University,}

\centerline{\it Beijing 100875, People's Republic of China}

\smallskip

\centerline{\it E-mail: skchen@mail.bnu.edu.cn and lizh@bnu.edu.cn}

\bigskip

\bigskip

\noindent{\bf Abstract.} For general (1+1)-affine Markov processes, we prove the ergodicity and
exponential ergodicity in total variation distances. Our methods follow the arguments of ergodic
properties for L\'{e}vy-driven OU-processes and a coupling of CBI-processes constructed by stochastic equations driven by time-space noises. Then the strong Feller property is considered.

\bigskip


\bigskip

\section{Introduction}

\noindent
Let $m\ge 0$ and $n\ge 0$ be integers. A time-homogeneous Markov processes $\{X_t:t\geq0\}=\{(Y_t,Z_t):t\geq0\}$ taking values in $\mbb{G}:=\mbb{R}_{+}^m\times\mbb{R}^n$ is called an affine Markov process if its characteristic function satisfies
 \beqlb\label{1.1}
\mbb{E}_x[\mrm{e}^{\langle X_t,u\rangle}]
 =
\exp\Big\{\langle x, V(t,u)\rangle+\int_0^t\psi(V(s,u))\,\mrm{d}s\Big\},
 \qquad
x\in\mbb{G}, u\in i\mbb{R}^{m+n},
 \eeqlb
where $V$ and $\psi$ are two complex-valued functions and $V$ satisfies certain generalized Riccati equations. The affine property means roughly that the logarithm of the characteristic function is affine with respect to the initial state. The concept of affine Markov processes enables a unified treatment of two important Markov classes including {\it continuous state branching processes with immigration} (CBI-processes) and {\it Ornstein-Uhlenbeck type processes} (OU-type processes), where the OU-type processes also include L\'evy processes as a particular case. Roughly speaking, affine processes with state space $\mbb{R}_{+}^m$ are $m$-dimensional CBI-processes, and those with state space $\mbb{R}^n$ are $n$-dimensional OU-type processes. The processes involve rich common mathematical structures and have found interesting connections and applications in several areas.
The general theory of finite-dimensional affine Markov processes including several equivalent characterizations and common financial applications was given by Duffie et al. (2003) under a regularity assumption, which requires the functions
 $$
t\mapsto V(t,u) \quad\mbox{and}\quad t\mapsto \int_0^t\psi(V(s,u))\,\mrm{d}s
 $$
are differentiable at $t=0$ and continuous at $u=0$. The regularity problem asks whether this property holds automatically for stochastically continuous affine processes. This property was established in Dawson and Li (2006) under the first moment condition. The problem was finally settled in Keller-Ressel et al. (2011), where it was proved that any stochastically continuous affine Markov process is regular. The connection of the regularity problem with Hilbert's fifth problem was also explained in Keller-Ressel et al. (2011).

The ergodicity and strong Feller property of CBI- and OU-type processes have been studied by a number of authors. In particular, a sufficient and necessary integrability condition for the ergodicity of a one-dimensional subcritical or critical CBI-process was announced in Pinsky (1972); see Li (2011) for a proof. It was proved in Sato and Yamazato (1984) that a finite-dimensional OU-type process is ergodic if and only if the eigenvalues of its coefficient matrix have strictly negative real parts. The coupling property and strong Feller property of finite-dimensional OU-type processes was studied in Priola and Zabczyk (2009) and Wang (2011). The ergodicity and exponential ergodicity of such processes in total variation distances were proved in Schilling and Wang (2012) and Wang (2012). The strong Feller property and exponential ergodicity in the total variation distance of one-dimensional CBI-processes were shown in Li and Ma (2015) by a coupling method; see also Li (2020a). In the recent work of Li (2020b), the ergodicities and exponential ergodicities in Wasserstein and total variation distances of Dawson-Watanabe superprocesses with or without immigration were proved, which clearly includes the finite-dimensional CBI-processes.

For general finite-dimensional affine Markov processes with the strictly negative real parts of eigenvalues of its coefficient matrix, a sufficient condition for ergodicity in weak convergence was given in Jin et al. (2020). The necessity of the condition was not established in Jin et al. (2020), so their result partially covers the those in Pinsky (1972), Li (2011) and Sato and Yamazato (1984). The exponential ergodicity of finite-dimensional affine Markov processes in the Wasserstein distance was established in Friesen et al. (2020). Zhang and Glynn (2018) provided  sufficient conditions for ergodicity and exponential ergodicity of such processes in the total variation distance.

The main purpose of this paper is to study the ergodicity of the finite-dimensional affine processes. For simplicity, we focus on the (1+1)-dimensional affine processes. We prove some results on the ergodicity and exponential ergodicity of the processes in total variation distances under natural conditions. Instead of Zhang and Glynn (2018), our approach is based on coupling methods developed in Schilling and Wang (2012) and Wang (2012); see also Li and Ma (2015), which answers the question appeared on Jin et al. (2017), pp.1145; see also Friesen and Jin (2020), pp.646.

The remainder of this paper is organized as follows. In section 2, we prove the ergodicity of (1+1)-dimensional affine processes in total variance distances. The exponential ergodic property in total variance distances is established under stronger conditions in section 3. Finally, the strong Feller property is studied in section 4.

\section{Ergodicity in total variance distances}

\setcounter{equation}{0}

In this and next section, we mainly prove the ergodicity and exponential ergodicity of (1+1)-affine Markov processes in total variance distances by dividing the progress of two processes issued from different points into two parts by a natural coupling of CBI-processes. It is well known that the coupling property along with the existence of a stationary measure can yield the ergodicity for the process; see, e.g., Wang (2012), which motivates the basic proof in this section, as well as for section 3. More precisely, we prove the coupling property for the process in Proposition \ref{t3.3} under proper conditions, and show the existence of a stationary distribution in Lemma \ref{t3.4}. Based on those results, we prove the ergodic property in Theorem \ref{t3.5}.

We first give some notations. Write $\mbb{G}=\mbb{R}_{+}\times\mbb{R}$.
Define
$U=\mbb{C}_{-}\times i\mbb{R}$, where $\mbb{C}_{-}=\{a+ib:a\in\mbb{R}_{-},b\in\mbb{R}\}$ and $i\mbb{R}=\{ia:a\in\mbb{R}\}.$ We further define two functions on $U$ as follows:
\beqnn
&&\phi(u)=-a_1u_1-b_1u_2+(\alpha_{11}+\alpha_{12})u_1^2+2(\sqrt{\alpha_{11}\alpha_{21}}+\sqrt{\alpha_{12}\alpha_{22}})u_1u_2\\
&&~~~~~~~~~~
+(\alpha_{21}+\alpha_{22})u_2^2+\int_{\mbb{G}}(\mrm{e}^{\langle u,z\rangle}-1-\langle
u,z\rangle)\,m(\mrm{d}z),\\
&&\psi(u)=a_2u_1-b_0u_2+\frac{1}{2}\sigma^2u_2^2
+\int_{\mbb{G}}(\mrm{e}^{\langle u,z\rangle}-1-z_2u_2)\,n(\mrm{d}z),
\eeqnn
where $a_1,b_i\in\mbb{R}$ $(i=0,1,2)$, $a_2, \sigma\in\mbb{R}_{+}$, $(\alpha_{ij})_{2\times2}$ is a nonnegative matrix, $m$ and $n$ are two L\'evy measures supported on $\mbb{G}\setminus\{0\}$ satisfying
 \beqlb
\int_{\mbb{G}}\Big(z_1\wedge z_1^2+|z_2|^2\Big)\,m(\mrm{d}z)
+\int_{\mbb{G}}\Big(1\wedge z_1+|z_2|\wedge|z_2|^2\Big)\,n(\mrm{d}z)<\infty.\label{3.1}
 \eeqlb
Denote by $\{P_t:t\geq0\}$ the transition semigroup of (1+1)-affine Markov process $\{X_t:t\geq0\}$. It is well known that $\{P_t:t\geq0\}$ can be uniquely determined by
$$
\int_{\mbb{G}}\mrm{e}^{\langle u,\xi\rangle}\,P_t(x,\mrm{d}\xi)=\exp\Big\{\langle x,
V(t,u)\rangle+\int_0^t\psi(V(s,u))\,\mrm{d}s\Big\},
$$
where
$$
\left\{\begin{array}{ll}
{\frac{\partial V_1}{\partial t}(t,u)
=\phi(V(t,u)),}
&V_1(0,u)=u_1. \\
&\\
{V_2(t,u)=\mrm{e}^{-b_2t}u_2}. &\quad
\end{array}\right.
$$
We can also obtain the process as a unique strong solution to a stochastic integral equation system. Let $(\Omega,\mathcal{F},\mathcal{F}_t,\mbb{P})$ be a filtered probability space satisfying the usual hypotheses. Denote by $W_i(\mrm{d}s,\mrm{d}u),$ $i=0,1,2$ the $(\mathcal{F}_t)$-Gaussian white noises on $(0,\infty)^2$ with intensity $\mrm{d}s\mrm{d}u$, $M(\mrm{d}s,\mrm{d}u,\mrm{d}z)$ be an $(\mathcal{F}_t)$-Poisson random measure on $(0,\infty)^2\times\mbb{G}$ with intensity $\mrm{d}s\mrm{d}um(\mrm{d}z)$ and $N(\mrm{d}s,\mrm{d}z)$ an $(\mathcal{F}_t)$-Poisson random measure on $(0,\infty)\times\mbb{G}$ with intensity $\mrm{d}sn(\mrm{d}z)$, the corresponding compensated measures are defined by $\tilde{M}(\mrm{d}s,\mrm{d}u,\mrm{d}z)$ and $\tilde{N}(\mrm{d}s,\mrm{d}z)$. Let $P_i^M(t)=\int_0^t\int_0^{\infty}\int_{\mbb{G}}\,z_i\,M(\mrm{d}s,\mrm{d}u,\mrm{d}z)$ and $P_i^N(t)=\int_0^t\int_{\mbb{G}}\,z_i\,N(\mrm{d}s,\mrm{d}z)$ for $i=1,2$. We assume those random elements are independent of each other. Let $Y_0,Z_0$ be $\mathcal{F}_0$-measurable random variables and $Y_0\geq0$. Let us consider the following stochastic integral equation system:
 \beqlb
&&Y_t=Y_0+\int_0^t(a_2-a_1Y_s)\,\mrm{d}s+\sqrt{2\alpha_{11}}\int_0^t\int_0^{Y_s}\,W_1(\mrm{d}s,\mrm{d}u)
\nonumber\\
&&~~~~~~~~~~
+\sqrt{2\alpha_{12}}\int_0^t\int_0^{Y_s}\,W_2(\mrm{d}s,\mrm{d}u)
+\int_0^t\int_{\mbb{G}}z_1\,N(\mrm{d}s,\mrm{d}z)\nonumber\\
&&~~~~~~~~~~
+\int_0^t\int_0^{Y_{s-}}\int_{\mbb{G}}z_1\,\tilde{M}(\mrm{d}s,\mrm{d}u,\mrm{d}z),\label{3.2}\\
&&Z_t=Z_0-\int_0^t(b_0+b_1Y_s+b_2Z_s)\,\mrm{d}s+\sigma\int_0^t\int_0^1\,W_0(\mrm{d}s,\mrm{d}u)\nonumber\\
&&~~~~~~~~~~
+\sqrt{2\alpha_{21}}\int_0^t\int_0^{Y_s}\,W_1(\mrm{d}s,\mrm{d}u)
+\sqrt{2\alpha_{22}}\int_0^t\int_0^{Y_s}\,W_2(\mrm{d}s,\mrm{d}u)\nonumber\\
&&~~~~~~~~~~
+\int_0^t\int_{\mbb{G}}z_2\,\tilde{N}(\mrm{d}s,\mrm{d}z)
+\int_0^t\int_0^{Y_{s-}}\int_{\mbb{G}}z_2\,\tilde{M}(\mrm{d}s,\mrm{d}u,\mrm{d}z).\label{3.3}
 \eeqlb
Here and in the sequel, we understand that for any $a\leq b\in\mbb{R}$
 $$
\int_a^b=\int_{(a,b]}\quad \mrm{and}
\quad \int_a^{\infty}=\int_{(a,\infty)}.
 $$
The existence of the solution to (\ref{3.2})--(\ref{3.3}) is a consequence of Theorem 6.2 in Dawson and Li (2006), where a weakly equivalent stochastic equation system was studied. The pathwise uniqueness for (\ref{3.2})--(\ref{3.3}) follows by modifications of the proofs in Dawson and Li (2006, 2012). Denote by $\{X_t:t\geq0\}= \{(Y_t,Z_t): t\geq0\}$ the unique strong solution to (\ref{3.2})--(\ref{3.3}). Then $\{Y_t:t\geq0\}$ is a one-dimensional CBI-process and $\{X_t:t\geq0\}$ is an $(1+1)$-dimensional affine Markov process.

For $x=(x_1,x_2)$ and $y=(y_1,y_2)\in\mbb{G}$, let $\{X_t(x):t\geq0\}$ and $\{X_t(y):t\geq0\}$ be the affine processes defined by (\ref{3.2})--(\ref{3.3}) starting from $x$ and $y$, respectively. Let $\varsigma=\inf\{t\geq0:Y_t(x_1)=Y_t(y_1)\}$ be the coalescence time of the coupling $\{(Y_t(x_1),Y_t(y_1)):t\geq0\}$. Given $f\in B_b(\mbb{G}),$ one can see that
 \beqlb
|P_tf(x)-P_tf(y)|
 &\leq&
\mbb{P}\Big\{|f(X_t(x))-f(X_t(y))|\mbf{1}_{\{t<\varsigma\}}\Big\}\nonumber\\
&&
+\Big|\mbb{P}\{[f(X_t(x))-f(X_t(y))]\mbf{1}_{\{t\geq\varsigma\}}\}\Big| \nonumber\\
 &\leq&
2\|f\|\mbb{P}\Big\{Y_t(x_1)\neq Y_t(y_1)\Big\}\nonumber\\
&&
+ \Big|\mbb{P}\{[f(X_t(x))-f(X_t(y))] \mbf{1}_{\{t\geq\varsigma\}}\}\Big|. \label{3.4}
 \eeqlb
For any $\varepsilon>0,$ we define a finite measure $n_{\varepsilon} $ on $\mbb{R}$ such that
 $$
n_{\varepsilon}(B)=\left\{\begin{array}{ll}
{n(\mbb{R}_{+}\times B)~,}
&\mbox{ if }\,n(\mbb{R}_{+}\times\mbb{R})<\infty; \\
&\\
{n(\mbb{R}_{+}\times B_{\varepsilon})~,} &\mbox{ if }\,n(\mbb{R}_{+}\times\mbb{R})=\infty,
\end{array}\right.
 $$
where $B\in\cal{B}(\mbb{R})$ and $B_{\varepsilon}:=B\setminus\{z_2:|z_2|<\varepsilon\}$.

In the following we give some key estimates. For $x=(x_1,x_2)$ and $y=(y_1,y_2)\in\mbb{G}$, let $\{X_t(x):t\geq0\}$ and $\{X_t(y):t\geq0\}$ be the affine processes defined by (\ref{3.2})--(\ref{3.3}) starting from $x$ and $y$, respectively. In view of (\ref{3.2})--(\ref{3.3}) we have
 \beqlb\label{3.5}
&&Z_t(x)=T_tx_2-\int_0^tT_{t-s}\big[b_0+b_1Y_s(x_1)\big]\,\mrm{d}s + \sqrt{2\alpha_{21}}\int_0^t\int_0^{Y_s(x_1)}T_{t-s}\,W_1(\mrm{d}s,\mrm{d}u)\nonumber\\
&&~~~~~~~~~~
+ \sqrt{2\alpha_{22}}\int_0^t\int_0^{Y_s(x_1)}T_{t-s}\,W_2(\mrm{d}s,\mrm{d}u) + \sigma\int_0^t\int_0^1T_{t-s}\,W_0(\mrm{d}s,\mrm{d}u) \nonumber\\
&&~~~~~~~~~~
+ \int_0^t\int_0^{Y_{s-}(x_1)}\int_{\mbb{G}}T_{t-s}z_2\,\tilde{M}(\mrm{d}s,\mrm{d}u,\mrm{d}z)
+ \int_0^t\int_{\mbb{G}}T_{t-s}z_2\,\tilde{N}(\mrm{d}s,\mrm{d}z),
 \eeqlb
and
 \beqlb\label{3.6}
&&Z_t(x)-Z_t(y)=T_t(x_2-y_2)-b_1\int_0^tT_{t-s}\Big(Y_s(x_1)-Y_s(y_1)\Big)\,\mrm{d}s\nonumber\\
&&~~~~~~~~~~~~~~~~~~~~~~~~
+\sqrt{2\alpha_{21}}\int_0^t\int_{Y_s(y_1)}^{Y_s(x_1)}T_{t-s}\,W_1(\mrm{d}s,\mrm{d}u)\nonumber\\
&&~~~~~~~~~~~~~~~~~~~~~~~~
+\sqrt{2\alpha_{22}}\int_0^t\int_{Y_s(y_1)}^{Y_s(x_1)}T_{t-s}\,W_2(\mrm{d}s,\mrm{d}u)\nonumber\\
&&~~~~~~~~~~~~~~~~~~~~~~~~
+\int_0^t\int_{Y_{s-}(y_1)}^{Y_{s-}(x_1)}\int_{\mbb{G}}T_{t-s}z_2\,\tilde{M}(\mrm{d}s,\mrm{d}u,\mrm{d}z),
 \eeqlb
where $T_t=\mrm{e}^{-b_2t}$ for $t\geq0$.

\bglemma\label{t3.1}
Suppose $0<2b_2<a_1$. Then there exist strictly positive constants $C_1$ and $C_2$ such that
for any $x,y\in\mbb{G}$
 \beqlb
&&\mbb{E}|Z_t(x)-Z_t(y)|\leq T_t\Big(|x_2-y_2|+C_2(x_1-y_1)+\sqrt{C_1(x_1-y_1)}\Big),\label{3.7}\\
&&\mbb{P}\{|Z_t(x)-Z_t(y)|>T_t\eta\}\leq\frac{1}{\eta}\Big(|x_2-y_2|+C_2(x_1-y_1) + \sqrt{C_1(x_1-y_1)}\Big)\label{3.8}
 \eeqlb
for $\eta>0$.
\edlemma

\proof
Note that $\mbb{E}[Y_s(x_1)-Y_s(y_1)]=(x_1-y_1)\mrm{e}^{-a_1s}.$ By Martingale inequality and Cauchy-Schwartz inequality we see that
 \beqnn
&&\mbb{E}\Big|\int_0^t\int_{Y_{s-}(y_1)}^{Y_{s-}(x_1)}
\int_{\mbb{G}}\mrm{e}^{b_2s}z_2\,\tilde{M}(\mrm{d}s,\mrm{d}u,\mrm{d}z)\Big|
\leq\bigg[\mbb{E}\Big|\int_0^t\int_{Y_{s-}(y_1)}^{Y_{s-}(x_1)}\int_{\mbb{G}}\mrm{e}^{b_2s}z_2\,
\tilde{M}(\mrm{d}s,\mrm{d}u,\mrm{d}z)\Big|^2\bigg]^{\frac{1}{2}}\\
&&~~~~~~~~~~~~~~~~~~~~~~~~~~~~~~~~~~~~~~~~~~~~~~~~~~~~~~~~~~~~~~~~~
\leq\bigg[\frac{\mrm{e}^{(2b_2-a_1)t}-1}{2b_2-a_1}(x_1-y_1)\int_{\mbb{G}}|z_2|^2\,m(\mrm{d}z)\bigg]^{\frac{1}{2}}\\
&&~~~~~~~~~~~~~~~~~~~~~~~~~~~~~~~~~~~~~~~~~~~~~~~~~~~~~~~~~~~~~~~~~
\leq\sqrt{C_{11}(x_1-y_1)},
 \eeqnn
similarly,
 \beqnn
&&\mbb{E}\Big|\sqrt{2\alpha_{21}}\int_0^t\int_{Y_s(y_1)}^{Y_s(x_1)}\mrm{e}^{b_2s}\,W_1(\mrm{d}s,\mrm{d}u)
+\sqrt{2\alpha_{22}}\int_0^t\int_{Y_s(y_1)}^{Y_s(x_1)}\mrm{e}^{b_2s}\,W_2(\mrm{d}s,\mrm{d}u)\Big|\nonumber\\
&&\leq\bigg[2\alpha_{21}\frac{\mrm{e}^{(2b_2-a_1)t}-1}{2b_2-a_1}(x_1-y_1)\bigg]^{\frac{1}{2}}+
\bigg[2\alpha_{22}\frac{\mrm{e}^{(2b_2-a_1)t}-1}{2b_2-a_1}(x_1-y_1)\bigg]^{\frac{1}{2}}\\
&&\leq\sqrt{C_{12}(x_1-y_1)},
 \eeqnn
and
 $$
\mbb{E}
\Big|-b_1\int_0^t\mrm{e}^{b_2s}(Y_{s-}(x_1)-Y_{s-}(y_1))\,\mrm{d}s\Big|
=\Big|\frac{b_1(\mrm{e}^{(b_2-a_1)t}-1)}{b_2-a_1}(x_1-y_1)\Big|\leq C_2(x_1-y_1),
 $$
where
$C_{11}=\frac{1}{a_1-2b_2}\int_{\mbb{G}}|z_2|^2\,m(\mrm{d}z), C_{12}=8\max\Big\{\frac{\alpha_{21}}{a_1-2b_2},\frac{\alpha_{22}}{a_1-2b_2}\Big\}$
and $C_2=\frac{|b_1|}{a_1-b_2}.$
By (\ref{3.6}) it is not hard to see that
 \beqnn
\mbb{E}|Z_t(x)-Z_t(y)|\leq T_t\Big(|x_2-y_2|+C_2(x_1-y_1)+\sqrt{C_1(x_1-y_1)}\Big),
 \eeqnn
where $C_1=4\max\{C_{11},C_{12}\}.$ The second inequality is an immediate result following from Markov inequality and (\ref{3.7}).\qed

For a bounded measurable function $f$ on $\mbb{R}$, define the supremum norm $\|f\|= \sup_x|f(x)|$. Given two bounded measures $\mu$ and $\nu$ on $(\mbb{R},\cal{B}(\mbb{R}))$, let $\mu\wedge\nu= \mu-(\mu-\nu)^{+}= \nu-(\nu-\mu)^{+}$, where the superscript ''$+$'' refers to the positive part in the Jordan-Hahn decomposition. It is easy to see that $\mu\wedge\nu= \nu\wedge\mu= 2^{-1}(\mu + \nu - |\mu-\nu|)$, where $|\mu-\nu|= (\mu-\nu)^{+} + (\nu-\mu)^{+}$ is the total variation measure. Let $\|\cdot\|_{\mrm{var}}$ denote the total variation norm defined by
  $$
\|\mu-\nu\|_{\mrm{var}}=\sup_{\|f\|\leq1}|\mu(f)-\nu(f)|,
  $$
where $\mu(f)= \int f \mrm{d}\mu$. For the convenience, we formulate the following conditions:

\noindent$ (\textbf{A})$ Denote the branching mechanism of the CBI-processes $\{Y_t: t\geq0\}$ by
$$
\phi_0(x)=a_1x+(\alpha_{11}+\alpha_{12})x^2+\int_{\mbb{G}}(\mrm{e}^{-xz_1}-1+xz_1)\,m(\mrm{d}z),\quad x\geq0.
$$
There exists $\theta>0$ such that for any $ z\geq\theta, \phi_0(z)>0$, and
 $$
 \int_{\theta}^{\infty}\phi_0^{-1}(z)\,\mrm{d}z<\infty.
 $$

\noindent$ (\textbf{B})$ There exist two constants $\varepsilon,~\eta>0$ such that
 \beqnn
\inf_{|a|\leq\eta}n_{\varepsilon}\wedge
(\delta_{a}\ast n_{\varepsilon}) (\mbb{R})>0,
\quad \int_{\{|z_2|>\varepsilon\}}|z_2|\,
n(\mrm{d}z)<\infty.
 \eeqnn

\begin{rem1}
\quad

(1). Condition (\textbf{A}) is called Grey's condition; see, e.g., Grey (1974). This condition has been used to study the exponentially ergodic property in total variance distances and strong Feller property of one-dimensional subcritical CBI-processes; see,e.g., Li and Ma (2015) and Li (2020a).

(2). Condition (\textbf{B}) is sharp to study the coupling property for the L\'evy process $(L_t)_{t\geq0}$ with L\'evy measure n and the Ornstein-Uhlenbeck process driven by $(L_t)_{t\geq0}$. Intuitively, it is one possibility to guarantee the sufficient jump activity such that the process admits a successful coupling; see,e.g., Schilling and Wang (2011,2012) and Wang (2012).
\end{rem1}

\bgproposition\label{t3.3}
Suppose that Conditions $\textbf{(A,B)}$ hold and $0<2b_2<a_1$. Then there exist constants $\hat{C},\kappa>0$, such that for any $t>0,x,y\in\mbb{G},$ we have
 $$
\|P_t(x,\cdot)-P_t(y,\cdot)\|_{\mrm{var}}\leq\hat{C}\Big(1+(\bar{v}_{\frac{1}
{\kappa^2+1}t}+1)|x_1-y_1|+\sqrt{|x_1-y_1|}+|x_2-y_2|\Big)\frac{1}{\sqrt{t}},
 $$
where $\bar{v}_t$ is the unique solution of the following equation:
 \beqlb
\frac{\mrm{d}}{\mrm{d}t}\bar{v}_t=-\phi_0(\bar{v}_t), \quad t>0 \label{3.9}
 \eeqlb
with initial condition $\bar{v}_{0+}=\infty$. Moreover, the mapping $t\mapsto\bar{v}_t$ is decreasing.
\edproposition

\proof
Under Condition \textbf{(A)}, the solution $\bar{v}_t$ to (\ref{3.9}) is unique, and the mapping $t\mapsto\bar{v}_t$ is decreasing; see, e.g., Theorem
3.6, Theorem 3.7 and Corollary 3.14 in Li (2020a). Under Condition \textbf{(B)}, for simplicity, we let
  $$
L_t^{\varepsilon}=\int_0^t\int_{\{|z_2|>\varepsilon\}}z_2\,N(\mrm{d}s,\mrm{d}z)
  $$
and define a sequence of stopping times $\Upsilon^{\varepsilon}_n= \inf\{t>\Upsilon^{\varepsilon}_{n-1}: L_t^{\varepsilon}\neq L_{t-}^{\varepsilon}\}$ with convention $\Upsilon^{\varepsilon}_0=0$. For $i\geq1$ let $\tau^{\varepsilon}_i= \Upsilon^{\varepsilon}_{i}-\Upsilon^{\varepsilon}_{i-1}$ and
$U^{\varepsilon}_i= \int_{\{\Upsilon^{\varepsilon}_{i}\}} \int_{\{|z_2|>\varepsilon\}} z_2\, N(\mrm{d}s,\mrm{d}z).$
Then $(\tau^{\varepsilon}_i)_{i\geq1}$ are i.i.d. random variables which are exponentially distributed with intensity $C_{\varepsilon}= n_{\varepsilon}(\mbb{R})$ and $(U^{\varepsilon}_i)_{i\geq1}$ of i.i.d. random variables on $\mbb{R}$ with distribution $\bar{n}_{\varepsilon}:=n_{\varepsilon}/C_{\varepsilon}$. Moreover, the two sequences $(U^{\varepsilon}_i)_{i\geq1}$ and $(\tau^{\varepsilon}_i)_{i\geq1}$ are
independent of each other. Let $N^{\varepsilon}_t=\sup\{k\geq1:\sum_{i=1}^k\tau^{\varepsilon}_i\leq t\}$. Then $(N^{\varepsilon}_t)_{t\geq0}$ is a Poisson process of intensity $C_{\varepsilon}$. Now we can rewrite
  $$
\int_0^t\int_{\{|z_2|>\varepsilon\}}z_2\,N(\mrm{d}s,\mrm{d}z)
 =
\sum\limits_{i=1}^{N^{\varepsilon}_t} U^{\varepsilon}_i
  $$
with $\sum_{i=1}^0=0$ by convention. It is not hard to check that
 \beqnn
\int_0^t\int_{\{|z_2|\geq\varepsilon\}}T_{t-s}z_2\,
N(\mrm{d}s,\mrm{d}z)=0\cdot\mbf{1}_{\{\Upsilon^{\varepsilon}_1>t\}}+\sum_{k=1}^{\infty}
\mbf{1}_{\{\Upsilon^{\varepsilon}_k\leq t<\Upsilon^{\varepsilon}_{k+1}\}}\sum_{i=1}^kT_{t-\Upsilon^{\varepsilon}_i}U^{\varepsilon}_i.
 \eeqnn

Let us use the following notations for the convenience,
 \beqnn
&&\zeta_{t,\varepsilon}=T_t
\Big\{ (b_0-\int_{\{|z_2|>\varepsilon\}}z_2\,n(\mrm{d}z))\int_0^t\mrm{e}^{b_2s}\,\mrm{d}s\\
&&~~~~~~~~~~~~~~~~
+\sigma\int_0^t\int_0^1\mrm{e}^{b_2s}\,
W_0(\mrm{d}s,\mrm{d}u)+\int_0^t\int_{\{|z_2|\le \varepsilon\}}\mrm{e}^{b_2s}z_2\,
\tilde{N}(\mrm{d}s,\mrm{d}z)\Big\},\\
&&\theta_t(x_1)=-b_1\int_0^t\mrm{e}^{b_2s}Y_s(x_1)\,
\mrm{d}s+
\int_0^t\int_0^{Y_{s-}(x_1)}\int_{\mbb{G}}\mrm{e}^{b_2s}z_2\,
\tilde{M}(\mrm{d}s,\mrm{d}u,\mrm{d}z),\\
&&~~~~~~~~~~~
+\sqrt{2\alpha_{21}}\int_0^t\int_0^{Y_s(x_1)}\mrm{e}^{b_2s}\,W_1(\mrm{d}s,\mrm{d}u)
+\sqrt{2\alpha_{22}}\int_0^t\int_0^{Y_s(x_1)}\mrm{e}^{b_2s}\,W_2(\mrm{d}s,\mrm{d}u).
 \eeqnn
Moreover, for any $t\geq0$, we denote the distribution of $\theta_t(x_1)-\theta_t(y_1)$ by $\Gamma_{t,x_1-y_1}$. In view of (\ref{3.3}), for any $f\in B_b(\mbb{G})$,
 \beqnn
&&\mbb{E}f(Y_t(x_1),Z_t(x))=\mbb{E}f\Big(Y_t(x_1),
T_t(x_2+\theta_t(x_1))+\zeta_{t,\varepsilon}+\sum^{N^{\varepsilon}_t}_{k=1}T_{t-\Upsilon^{\varepsilon}_k}U^{\varepsilon}_k
\Big)\\
&&=\mbb{E}\Big[f\Big(Y_t(x_1),T_t(x_2+\theta_t(x_1))+\zeta_{t,\varepsilon}
\Big)\mbf{1}_{\{N^{\varepsilon}_t=0\}}\Big]\\
&&~~
+\mbb{E}f\Big(Y_t(x_1),\sum_{k=1}^{\infty}\mbf{1}_{\{\Upsilon^{\varepsilon}_k\leq t<\Upsilon^{\varepsilon}_{k+1}\}}
\Big\{T_t(x_2+\theta_t(x_1))
+\zeta_{t,\varepsilon}+\sum_{j=1}^{k}T_{\Upsilon^{\varepsilon}_j}U^{\varepsilon}_j\Big\}\Big)\\
&&=\mbb{E}\Big[f\Big(Y_t(x_1),T_t(x_2+\theta_t(x_1))+\zeta_{t,\varepsilon}\Big)\mbf{1}_{\{N^{\varepsilon}_t=0\}}\Big]\\
&&~~+\sum_{k=1}^{\infty}\int ...\int_{\sum_{i=1}^{k}t_i\leq t<\sum_{i=1}^{k+1}t_i}C^{k+1}_{\varepsilon}
\mrm{e}^{-C_{\varepsilon}\sum_{i=1}^{k+1}t_i}\,\mrm{d}t_1\,...\mrm{d}t_{k+1}\\
&&~~
\times\int_{\mbb{R}^k}\mbb{E}f\Big(Y_t(x_1),
T_t(x_2+\theta_t(x_1))+\zeta_{t,\varepsilon}+\sum_{i=1}^{k}T_{\sum_{j=1}^it_j}r_i\Big)\,
\bar{n}_{\varepsilon}(\mrm{d}r_1)
\,...\bar{n}_{\varepsilon}(\mrm{d}r_k)\\
&&=\mbb{E}\Big[f\Big(Y_t(x_1),T_t(x_2+\theta_t(x_1))+\zeta_{t,\varepsilon}
\Big)\mbf{1}_{\{N^{\varepsilon}_t=0\}}\Big]\\
&&~~+\sum_{k=1}^{\infty}\int ...\int_{\sum_{i=1}^kt_i\leq t<\sum_{i=1}^{k+1}t_i}
C^{k+1}_{\varepsilon}\mrm{e}^{-C_{\varepsilon}\sum_{i=1}^{k+1}t_i}\,\mrm{d}t_1\,...\mrm{d}t_{k+1}\\
&&~~\quad
\times\int_{\mbb{R}} \mbb{E}f\Big(Y_t(x_1),
T_t(x_2+\theta_t(x_1))+\zeta_{t,\varepsilon}+z\Big)\,n_{t_1,...,t_k}(\mrm{d}z),
 \eeqnn
where the second equality partly follows from formula (2.10) in Schilling and Wang (2012); see also Lemma 2.2 in Schilling and Wang (2012). Here, $n_{t_1,...t_k}(\mrm{d}z)$ is the probability measure on $\mbb{R}$, which
is the image of the $k$-fold product measure
$\bar{n}_{\varepsilon}\times...\times\bar{n}_{\varepsilon}$
under the linear transformation $J_{t_1,...,t_k}:$
$J_{t_1,...t_k}(r_1,...,r_k)=T_{t_1}r_1+...+T_{t_1+...+t_k}r_k.$
Given $(x_1,x_2),(y_1,y_2)\in\mbb{G},$ without loss of generality, we can assume that $x_1\geq y_1$, and
we have
 \beqnn
&&\sup_{\|f\|\leq1}\Big|\mbb{E}\big[(f(X_t(x))-f(X_t(y)))
\mbf{1}_{\{t\geq\varsigma\}}\big]\Big|\\
&&\leq 2\mrm{e}^{-C_{\varepsilon}t}+\sum_{k=1}^{\infty}\int ...\int_{\sum_{i=1}^kt_i\leq t<\sum_{i=1}^{k+1}t_i}
C^{k+1}_{\varepsilon}\mrm{e}^{-C_{\varepsilon}\sum_{i=1}^{k+1}t_i}\,\mrm{d}t_1\,...\mrm{d}t_{k+1}\\
&&~~~~~~~~~~~~~~
\times\sup_{\|f\|\leq1}\Big|\int_{\mbb{R}}\mbb{E}f\Big(Y_t(x_1),
T_t(x_2+\theta_t(x_1))+\zeta_{t,\varepsilon}+z\Big)\,n_{t_1,...,t_k}(\mrm{d}z) \\
&&~~~~~~~~~~~~~~~~~~~
 - \int_{\mbb{R}}\mbb{E}f\Big(Y_t(y_1),
T_t(y_2+\theta_t(y_1))+\zeta_{t,\varepsilon}+z\Big)\,n_{t_1,...,t_k}(\mrm{d}z)\Big|\\
&&\leq2\mrm{e}^{-C_{\varepsilon}t}+\sum_{k=1}^{\infty}\int ...\int_{\sum_{i=1}^kt_i\leq t<\sum_{i=1}^{k+1}t_i}
C^{k+1}_{\varepsilon}\mrm{e}^{-C_{\varepsilon}\sum_{i=1}^{k+1}t_i}\,\mrm{d}t_1\,...\mrm{d}t_{k+1}\\
&&~~~~~~~~~~~~~~
\times\sup_{\|f\|\leq1}\int_{\mbb{R}}\,\Big|\int_{\mbb{R}}\mbb{E}f(Y_t(y_1),
T_t(y_2+\theta_t(y_1))+\zeta_{t,\varepsilon}+z)\, \\
&&~~~~~~~~~~~~~~~~~~~~~~~~~~~~~~~~~\qqquad
\times\delta_{T_t(x_2-y_2+z_2)}\ast n_{t_1,...,t_k}(\mrm{d}z) \\
&&~~~~~~~~~~~~~~~~~~
 - \int_{\mbb{R}}\mbb{E}f(Y_t(y_1),
T_t(y_2+\theta_t(y_1))+\zeta_{t,\varepsilon}+z\Big)\,n_{t_1,...,t_k}(\mrm{d}z)\Big|\,
\Gamma_{t,x_1-y_1}(\mrm{d}z_2)\\
&&\leq2\mrm{e}^{-C_{\varepsilon}t}+\sum_{k=1}^{\infty}\int ...\int_{\sum_{i=1}^kt_i\leq t<\sum_{i=1}^{k+1}t_i}
C^{k+1}_{\varepsilon}\mrm{e}^{-C_{\varepsilon}\sum_{i=1}^{k+1}t_i}\,\mrm{d}t_1\,...\mrm{d}t_{k+1}\\
&&~~~~~~~~~~~~~~
\times\int_{W}\varpi(z_2,t,t_1,...,t_k)\,\Gamma_{t,x_1-y_1}(\mrm{d}z_2)\\
&&~~~~~~~~~~~~
+2\sum_{k=1}^{\infty}\int ...\int_{\sum_{i=1}^kt_i\leq t<\sum_{i=1}^{k+1}t_i}C^{k+1}_{\varepsilon}
\mrm{e}^{-C_{\varepsilon}\sum_{i=1}^{k+1}t_i}\,\mrm{d}t_1\,...\mrm{d}t_{k+1}
\int_{\mbb{R}\setminus W}\,\Gamma_{t,x_1-y_1}(\mrm{d}z_2),
 \eeqnn
where $W=\{z_2:|z_2+x_2-y_2|< \eta\}$ for some $\eta>0$ and
 \beqnn
&&\varpi(z_2,t,t_1,...,t_k)=\|\delta_{T_t(x_2-y_2+z_2)}\ast
n_{t_1,...,t_k}-n_{t_1,...,t_k}\|_{\mrm{var}}.
 \eeqnn
If we set $p(t,x,y):=\int_{\mbb{R}\setminus W}\,\Gamma_{t,x_1-y_1}(\mrm{d}z_2),$ one can see that
 \beqnn
&&\sum_{k=1}^{\infty}\int...\int_{\sum_{i=1}^kt_i\leq t<\sum_{i=1}^{k+1}t_i}
C^{k+1}_{\varepsilon}\mrm{e}^{-C_{\varepsilon}\sum_{i=1}^{k+1}t_i}\,\mrm{d}t_1\,...\mrm{d}t_{k+1}
\int_{\mbb{R}\setminus W}\,
\Gamma_{t,x_1-y_1}(\mrm{d}z_2)\nonumber\\
&&\leq p(t,x,y)\mrm{e}^{-C_{\varepsilon}t}\sum_{k=1}^{\infty}C^{k+1}_{\varepsilon}
\int...\int_{t\geq\sum_{i=1}^kt_i}\,\mrm{d}t_1\,...\mrm{d}t_k\\
&&=p(t,x,y)C_{\varepsilon}(1-\mrm{e}^{-C_{\varepsilon}t}).
 \eeqnn
By (\ref{3.6}), we have$ T_t^{-1}[Z_t(x)-Z_t(y)]=(x_2-y_2) + \theta_t(x_1)-\theta_t(y_1).$
It follows that $T_t^{-1}|Z_t(x)-Z_t(y)|=|x_2-y_2+\theta_t(x_1)-\theta_t(y_1))|$. Then
 \beqlb\label{3.10}
&&p(t,x,y)\le\mbb{P}\{|Z_t(x)-Z_t(y)|>T_t\eta\}\nonumber\\
&&~~~~~~~~~~~~~~
\leq\frac{1}{\eta}\Big(|x_2-y_2|+C_2(x_1-y_1) + \sqrt{C_1(x_1-y_1)}\Big),
 \eeqlb
where the last inequality follows from ({\ref{3.8}}). On the other hand,
following the proof of Theorem 1.1 in Schilling and Wang (2012), for
$|x_2-y_2+z_2|\leq\eta,t\geq t_1+...+t_k$ and $k\geq1$ we can find some constants $C_{1,\eta}, C_3>0$ such that
 $$
\varpi(z_2,t,t_1,...,t_k)
\leq\frac{C_{1,\eta}}{\sqrt{k}},
 $$
and
 \beqnn
&&\frac{C_3}{\sqrt{t}}\geq2\mrm{e}^{-C_{\varepsilon}t}+\sum_{k=1}^{\infty}
\int ...\int_{\sum_{i=1}^kt_i\leq t<\sum_{i=1}^{k+1}t_i}C^{k+1}_{\varepsilon}
\mrm{e}^{-C_{\varepsilon}\sum_{i=1}^{k+1}t_i}\,\mrm{d}t_1\,...\mrm{d}t_{k+1}\\
&&~~~~~~~~~~~~~~~~
\times\int_{W}\varpi(z_2,t,t_1,...,t_k)\,\Gamma_{t,x_1-y_1}(\mrm{d}z_2).
 \eeqnn
In conclusion, we have
 \beqnn
&&\Big|\mbb{E}\big[(f(X_t(x))-f(X_t(y)))\mbf{1}_{\{t\geq\varsigma\}}\big]\Big|
\leq\frac{\|f\|C_3}{\sqrt{t}}+\frac{2\|f\|C_{\varepsilon}}{\eta}\\
&&~~~~~~~~~~~~~~~~~~~~\qqquad\qqquad\quad\times\Big(|x_2-y_2|+C_2(x_1-y_1)
+\sqrt{C_1(x_1-y_1)}\Big).
 \eeqnn
In view of (\ref{3.4}), for $f\in B_b(\mbb{G})$
 \beqnn
&&|P_tf(x)-P_tf(y)|\leq2\|f\|\bar{v}_t(x_1-y_1)+\frac{\|f\|C_3}{\sqrt{t}}\\
&&~~~~~~~~~~~~~~~~~~~~~~~~~~~~~~~~
+\frac{2\|f\|C_{\varepsilon}}{\eta}\Big(|x_2-y_2|+C_2(x_1-y_1)+\sqrt{C_1(x_1-y_1)}\Big).
 \eeqnn
For any $s\in(0,t),$ by using the Markov property and (\ref{3.7}),
 \beqnn
&&|P_tf(x)-P_tf(y)|\\
&&=\Big|\mbb{E}\Big\{P_sf(Y_{t-s}(x_1),Z_{t-s}(x))
-P_sf(Y_{t-s}(y_1),Z_{t-s}(y))\Big\}\Big|\\
&&\leq2\|f\|\bar{v}_s\mbb{E}(Y_{t-s}(x_1)-Y_{t-s}(y_1))
+\frac{\|f\|C_3}{\sqrt{s}}\\
&&~~~~
+\frac{2\|f\|C_{\varepsilon}}{\eta}\Big\{\mbb{E}|Z_{t-s}(x)-Z_{t-s}(y)|
+\Big[C_1\mbb{E}\{Y_{t-s}(x_1)-Y_{t-s}(y_1)\}\Big]^{\frac{1}{2}}\\
&&~~~~~~~~~~~~~~~~~~
+C_2\mbb{E}\{Y_{t-s}(x_1)-Y_{t-s}(y_1)\}\Big\}\\
&&\leq 2\|f\|\bar{v}_s(x_1-y_1)\mrm{e}^{-a_1(t-s)}+\frac{\|f\|C_3}{\sqrt{s}}\\
&&~~~~
+\frac{2\|f\|C_{\varepsilon}}{\eta}\Big\{\sqrt{C_1(x_1-y_1)\mrm{e}^{-a_1(t-s)}}+C_2(x_1-y_1)\mrm{e}^{-a_1(t-s)}\\
&&~~~~~~~~~~~~~~~~~~
+
\mrm{e}^{-b_2(t-s)}(|x_2-y_2|+C_2(x_1-y_1)+\sqrt{C_1(x_1-y_1)})\Big\},
 \eeqnn
which implies that
 \beqnn
&&\|P_t(x,\cdot)-P_t(y,\cdot)\|_{\mrm{var}}\\
&\leq& C_4\Big[1+(\bar{v}_s+1)(x_1-y_1)+\sqrt{x_1-y_1}+|x_2-y_2|\Big]
\Big(\mrm{e}^{-b_2(t-s)}\vee\frac{1}{\sqrt{s}}\Big)\\
&\leq& C_4\Big[1+(\bar{v}_s+1)(x_1-y_1)+\sqrt{x_1-y_1}+|x_2-y_2|\Big]
\Big(\frac{\kappa}{\sqrt{t-s}}\vee\frac{1}{\sqrt{s}}\Big)
 \eeqnn
for some constants $C_4>0$~and $\kappa=1\vee\frac{1}{\mrm{e}b_2}.$ We finally obtain the required assertion by setting $s=\frac{1}{\kappa^2+1}t,\,\,\hat{C}=\sqrt{\kappa^2+1}C_4.$
\qed

Note that the initial assumptions of L\'{e}vy measures $m$ and $n$ in Jin et al. (2020) are weaker than (\ref{3.1}). Before we establish the main result in this section, we need the following lemma, which is a consequence of Theorem 2.7 in Jin et al. (2020).

\bglemma\label{t3.4}
Assume that $a_1>0, b_2>0$, (\ref{3.1}) and
 \beqlb
 \int_{\{z_1\geq1\}}\log z_1\,n(\mrm{d}z)<\infty \label{3.11}
 \eeqlb
hold. Then the law of $X_t$ converges weakly to a limiting distribution $\pi$ given by
 $$
\int_{\mbb{G}} \mrm{e}^{\langle u,y\rangle}\,\pi(\mrm{d}y)
=\exp\Big\{\int_0^{\infty} \psi(V(s,u))\,\mrm{d}s\Big\},
\quad u\in U,
 $$
Moreover, $\pi$ is the unique stationary distribution for $X.$
\edlemma

Based on the Lemma \ref{t3.4} and Proposition \ref{t3.3} above, we can prove the ergodicity of the transition semigroup $(P_t)_{t\geq0}.$ More explicitly, we have the following result:

\bgtheorem\label{t3.5}
Suppose $\int_{\{z_1\geq1\}}\log z_1\,n(\mrm{d}z)< \infty$ and the conditions of Proposition \ref{t3.3} are satisfied. Then the affine process $(X_t)_{t\geq0}$ is ergodic in the total variation distance. Namely, there exists a unique invariant measure $\pi$ for the process such that for any $x\in\mbb{G},$
 \beqlb
\lim\limits_{t\to \infty}
\|P_t(x,\cdot)-\pi(\cdot)\|_{\mrm{var}}=0.\nonumber
 \eeqlb
\edtheorem

\proof
By (\ref{3.1}) and Lemma \ref{t3.4}, there exists a unique invariant measure $\pi$ for $(P_t)_{t\geq0}.$ Fix $x\in\mbb{G},$ one can see that
 $$
\|P_t(x,\cdot)-\pi(\cdot)\|_{\mrm{var}}\leq\int_{\mbb{G}}\|P_t(x,\cdot)-P_t(y,\cdot)\|_{\mrm{var}}\,
\pi(\mrm{d}y).
 $$
For any $\epsilon>0,$ we choose $\delta>0$ such that $\pi\{|x_1-y_1|+|x_2-y_2|>\delta\}\leq\epsilon.$
By Proposition \ref{t3.3} we have
 \beqnn
&&\|P_t(x,\cdot)-\pi(\cdot)\|_{\mrm{var}}\\
&\leq&\int_{\mbb{G}\setminus\Xi}\|P_t(x,\cdot)-P_t(y,\cdot)\|_{\mrm{var}}\,\pi(\mrm{d}y)
+\int_{\Xi}\|P_t(x,\cdot)-P_t(y,\cdot)\|_{\mrm{var}}\,\pi(\mrm{d}y)\\
&\leq&\frac{\hat{C}}{\sqrt{t}}\int_{\mbb{G}\setminus\Xi}
\Big(1+(\bar{v}_{\frac{1}{\kappa^2+1}t}+1)\delta+\sqrt{\delta}+\delta\Big)\,\pi(\mrm{d}y)+2\pi(\Xi)\\
&\leq&\frac{\hat{C}}{\sqrt{t}}\int_{\mbb{G}\setminus\Xi}
\Big(1+(\bar{v}_{\frac{1}{\kappa^2+1}t}+1)\delta+\sqrt{\delta}+\delta\Big)\,\pi(\mrm{d}y)+2\epsilon,
 \eeqnn
where $\Xi=\{y:|x_1-y_1|+|x_2-y_2|>\delta\}.$ Letting $t\to \infty$ first and then $\epsilon\to 0,$ we prove the desired result.\qed

\section{Exponential ergodicity in total variance distances}

\setcounter{equation}{0}

In this section, we prove the exponential ergodicity of (1+1)-affine Markov processes defined by (\ref{3.2})--(\ref{3.3}) with (\ref{3.1}) in total variance distance under proper conditions. Since the exponential ergodicity implies the ergodicity, it's reasonable to strengthen the assumptions to guarantee the exponential ergodicity. We also give the following conditions before moving forward:

\noindent$ (\textbf{C})$ There exists $\varepsilon>0,$ such that
 \beqnn
\limsup\limits_{\rho\to 0}
\Big[\frac{\sup_{|a|\leq\rho}\|n_{\varepsilon}
-\delta_{a}\ast n_{\varepsilon}\|_{\mrm{var}}}
{\rho}\Big]<\infty,
\quad \int_{\{|z_2|>\varepsilon\}}|z_2|\,
n(\mrm{d}z)<\infty.
 \eeqnn

\medskip

\noindent$ (\textbf{C}')$ There exists $\varepsilon>0,$ such that
 \beqnn
\limsup\limits_{\rho\to 0}
\Big[\frac{\sup_{|a|\leq\rho}\|n_{\varepsilon}
-\delta_{a}\ast n_{\varepsilon}\|_{\mrm{var}}}
{\rho}\Big]<\infty,
\quad \int_{\{|z_2|>1\}}|z_2|^2\,
n(\mrm{d}z)<\infty. \nonumber
 \eeqnn

\begin{rem1}
Condition (\textbf{C}') plays an important role in characterizing the exponential ergodicity. This condition implies $\int_{\mbb{G}}|z_2|^2\,n(\mrm{d}z)<\infty$, so it is stronger than (\textbf{C}). Condition (\textbf{B}) is weaker than (\textbf{C}); see, e.g., Remark~1 in Wang (2012) for a proof.

\end{rem1}

\bgproposition\label{t4.2} Suppose that Conditions $\textbf{(A,C)}$ hold and $0<2b_2<a_1$. Then there exist constants
$\tilde{C},\tilde{\kappa}>0,$ such that for any $t>0,x,y\in\mbb{G},$ we have
 $$
\|P_t(x,\cdot)-P_t(y,\cdot)\|_{\mrm{var}}\leq\tilde{C}
\Big(1+(\bar{v}_{\frac{\tilde{\kappa} t}{C_{\varepsilon}}}+1)|x_1-y_1|+\sqrt{|x_1-y_1|}+|x_2-y_2|\Big)
\mrm{e}^{-\tilde{\kappa} t},
 $$
where $\bar{v}_t$ is defined as that in Proposition \ref{t3.3}.
\edproposition

\proof
Following the arguments in the proof of Proposition \ref{t3.3}, it implies that on $\{N^{\varepsilon}_t\geq1\}$
 $$
\int_0^t\int_{\{|z_2|\geq\varepsilon\}}T_{t-s}z_2\,N(\mrm{d}s,\mrm{d}z)
=\sum_{k=1}^{N^{\varepsilon}_t}T_{t-\Upsilon^{\varepsilon}_k}U^{\varepsilon}_k.
 $$
For given $f\in B_b(\mbb{G}), x\in\mbb{G}$ we have the following decomposition
 $$
\mbb{E}[f\left(Y_t(x_1),Z_t(x)\right)]=\mbb{E}[f(Y_t(x_1),Z_t(x))\mbf{1}_{\{N^{\varepsilon}_t=0\}}]+P^1_tf(x),
 $$
where $P^1_tf(x)=\mbb{E}[f(Y_t(x_1),Z_t(x))\mbf{1}_{\{N^{\varepsilon}_t\geq1\}}]$ and
 \beqnn
&&P^1_tf(x)=\mbb{E}\Big\{\mbf{1}_{\{N^{\varepsilon}_t\geq1\}}f\Big(Y_t(x_1),
T_t(x_2+\theta_t(x_1))
+\zeta^1_{t,\varepsilon}+T_{t-\Upsilon^{\varepsilon}_{N^{\varepsilon}_t}}U^{\varepsilon}_{N^{\varepsilon}_t}\Big)\Big\}\\
&&~~~~~~~~~~~~
=\frac{1}{C_{\varepsilon}}\mbb{E}\Big\{\mbf{1}_{\{N^{\varepsilon}_t\geq1\}}\int_{\mbb{R}}f\Big(Y_t(x_1),
T_t(x_2+\theta_t(x_1))+\zeta^1_{t,\varepsilon}+T_{t-\Upsilon^{\varepsilon}_{N^{\varepsilon}_t}}z\Big)\,n_{\varepsilon}(\mrm{d}z)\Big\},
 \eeqnn
where
 \beqnn
&&\zeta^1_{t,\varepsilon}=T_t\Big\{(b_0-\int_{\{|z_2|>\varepsilon\}}z_2\,n(\mrm{d}z))
\int_0^t\mrm{e}^{b_2s}\,\mrm{d}s
+\sigma\int_0^t\int_0^1\mrm{e}^{b_2s}\,W_0(\mrm{d}s,\mrm{d}u)\\
&&~~~~~~~~~~~~~~~~~~
+\sum_{k=1}^{N^{\varepsilon}_t-1}T_{-\Upsilon^{\varepsilon}_k}U^{\varepsilon}_k
+\int_0^t\int_{\{|z_2|\leq\varepsilon\}}\mrm{e}^{b_2s}z_2\,\tilde{N}(\mrm{d}s,\mrm{d}z)\Big\}.
 \eeqnn
We recall again that the distributions of
$\theta_t(x_1)-\theta_t(y_1)$ is $\Gamma_{t,x_1-y_1}$. Given $x, y\in\mbb{G}$, without loss of generality, it suffices to
consider the case of $x_1\geq y_1.$
 \beqnn
&&\Big|\mbb{E}[(f(X_t(x)-f(X_t(y))))\mbf{1}_{\{N^{\varepsilon}_t\geq1\}}\mbf{1}_{\{t\geq\varsigma\}}]\Big|\\
&&\leq\frac{1}{C_{\varepsilon}}\Big|\mbb{E}\Big\{\mbf{1}_{\{N^{\varepsilon}_t\geq1\}}
\Big(\int_{\mbb{R}}f(z_1,t,y_1,y_2)\,n_{\varepsilon}\Big(\mrm{d}z_1
-T_{\Upsilon^{\varepsilon}_{N^{\varepsilon}_t}}(x_2-y_2+z_2)\Big)\\
&&~~~~~~
\times\,\int_{\mbb{R}}\,\Gamma_{t,x_1-y_1}(\mrm{d}z_2)-\int_{\mbb{R}}f(z_1,t,y_1,y_2)\,
n_{\varepsilon}(\mrm{d}z_1)\Big)\Big\}\Big|\\
&&\leq\frac{(1-\mrm{e}^{-C_{\varepsilon}t})\Lambda\,\|f\|}{C_{\varepsilon}}
\int_{\mbb{R}}\Big(|x_2-y_2|+|z_2|\Big)\,\Gamma_{t,x_1-y_1}(\mrm{d}z_2),
 \eeqnn
where
 \beqnn
&&\Lambda=\sup\limits_{\rho>0}\Big[\frac{\sup_{|a|\leq\rho}\|n_{\varepsilon}-\delta_{a}\ast n_{\varepsilon}\|_{\mrm{var}}}{\rho}\Big],\\
&&f(z_1,t,y_1,y_2)=f\Big(Y_t(y_1),T_t(y_2+\theta_t(y_1))+\zeta^1_{t,\varepsilon}+T_{t-\Upsilon^{\varepsilon}_{N^{\varepsilon}_t}}z_1\Big).
 \eeqnn
We have $\Lambda<\infty$ due to the Condition $\textbf{(C)}$ and the fact that $\sup\limits_{|a|\leq\rho}\|n_{\varepsilon}-\delta_{a}\ast n_{\varepsilon}\|_{\mrm{var}}\leq2C_{\varepsilon}.$ It together with Lemma \ref{t3.1} follow that
 $$
\Big|\mbb{E}[(f(X_t(x))-f(X_t(y)))\mbf{1}_{\{N^{\varepsilon}_t\geq1\}}\mbf{1}_{\{t\geq\varsigma\}}]\Big|
\leq\frac{\Lambda\|f\|}{C_{\varepsilon}}\Big\{|x_2-y_2|+\sqrt{C_1(x_1-y_1)}+C_2(x_1-y_1)\Big\}.
 $$
On the other hand,
 $$
\mbb{E}\Big|f(Y_t(x_1), Z_t(x))\mbf{1}_{\{N^{\varepsilon}_t=0\}}\mbf{1}_{\{t\geq\varsigma\}}\Big|
\leq\|f\|\mrm{e}^{-C_{\varepsilon}t},\quad t\geq0,f\in B_b(\mbb{G}),
 $$
it follows that
 \beqnn
&&\Big|\mbb{E}[(f(X_t(x))-f(X_t(y)))\mbf{1}_{\{t\geq\varsigma\}}]\Big|\\
&&\leq2\|f\|\mrm{e}^{-C_{\varepsilon}t}+\frac{\Lambda}{C_{\varepsilon}}\,\|f\|
\Big\{|x_2-y_2|+\sqrt{C_1(x_1-y_1)}+C_2(x_1-y_1)\Big\}.
 \eeqnn
Therefore, it together with (\ref{3.4}) and Theorem 10.3 in Li (2020a) imply that
 \beqnn
&&|P_tf(x)-P_tf(y)|\leq2\|f\|(\bar{v}_t(x_1-y_1)+\mrm{e}^{-C_{\varepsilon}t})\\
&&~~~~~~~~~~~~~~~~~~~~~~~~~~~~~~
+\|f\|\Big\{\frac{\Lambda}{C_{\varepsilon}}(|x_2-y_2|+\sqrt{C_1(x_1-y_1)}+C_2(x_1-y_1))\Big\}.
 \eeqnn
Similarly as that in Proposition \ref{t3.3}, for any $s\in(0,t),$ by using the Markov property and (\ref{3.7}) we have
 \beqnn
&&|P_tf(x)-P_tf(y)|\\
&&=\Big|\mbb{E}\Big[P_sf(Y_{t-s}(x_1),Z_{t-s}(x))-P_sf(Y_{t-s}(y_1),Z_{t-s}(y))\Big]\Big|\\
&&\leq2\|f\|\bar{v}_s\mbb{E}[Y_{t-s}(x_1)-Y_{t-s}(y_1)]+2\|f\|
\mrm{e}^{-C_{\varepsilon}s}\\
&&~~~~~~~
+\frac{\|f\|\Lambda}{C_{\varepsilon}}\Big\{
\mbb{E}|Z_{t-s}(x)-Z_{t-s}(y)|+\Big\{C_1\mbb{E}[Y_{t-s}(x_1)-Y_{t-s}(y_1)]\Big\}^{\frac{1}{2}}\\
&&~~~~~~~~~~~~~~~~~~
+C_2\mbb{E}[Y_{t-s}(x_1)-Y_{t-s}(y_1)]\Big\}\\
&&
\leq2\|f\|\bar{v}_s(x_1-y_1)\mrm{e}^{-a_1(t-s)}+2\|f\|
\mrm{e}^{-C_{\varepsilon}s}\\
&&~~~~~~~~~~~~
+\frac{\|f\|\Lambda}{C_{\varepsilon}}
\Big\{\mrm{e}^{-b_2(t-s)}\Big(|x_2-y_2|+C_2(x_1-y_1)\\\
&&~~~~~~~~~~~~
+\sqrt{C_1(x_1-y_1)}\Big)+[C_1(x_1-y_1)\mrm{e}^{-a_1(t-s)}]^{\frac{1}{2}}+C_2(x_1-y_1)\mrm{e}^{-a_1(t-s)}\Big\},
 \eeqnn
it follows that
 $$
\|P_t(x,\cdot)-P_t(y,\cdot)\|_{\mrm{var}}\leq\tilde{C}\Big(1+(\bar{v}_s+1)(x_1-y_1)+\sqrt{x_1-y_1}+|x_2-y_2|\Big)
\Big\{\mrm{e}^{-C_{\varepsilon}s}\vee \mrm{e}^{-b_2(t-s)}\Big\},
 $$
where
$\tilde{C}=\max\{2,\frac{\Lambda}{C_{\varepsilon}},\frac{C_2\Lambda}{C_{\varepsilon}},\frac{\sqrt{C_1}\Lambda}{C_{\varepsilon}}\}.$
Setting $s=\frac{b_2t}{C_{\varepsilon}+b_2},$ we obtain
 $$
\|P_t(x,\cdot)-P_t(y,\cdot)\|_{\mrm{var}}
\leq\tilde{C}\Big(1+(\bar{v}_{\frac{\tilde{\kappa} t}
{C_{\varepsilon}}}+1)(x_1-y_1)+\sqrt{x_1-y_1}+|x_2-y_2|\Big)\mrm{e}^{-\tilde{\kappa} t},
 $$
where $\tilde{\kappa}=\frac{b_2C_{\varepsilon}}{C_{\varepsilon}+b_2}.$ Then the required assertion holds.
\qed

Based on Proposition \ref{t4.2}, we have the following result:

\bgtheorem\label{t4.3} Suppose that Conditions $\textbf{(A,C')}$ hold and $0<2b_2<a_1$. Moreover,
assume that $\int_{\{z_1>1\}}z_1\,n(\mrm{d}z)<\infty$ holds. Then the affine
process $(X_t)_{t\geq0}$ is exponentially ergodic in the total variation distance.
\edtheorem

\proof
Given $x\in\mbb{G},$ it follows from Proposition \ref{t4.2} that
 \beqnn
&&\|P_t(x,\cdot)-\pi(\cdot)\|_{\mrm{var}}
\leq\int_M\|P_t(x,\cdot)-P_t(y,\cdot)\|_{\mrm{var}}\,\pi(\mrm{d}y)\\
&&~~~~~~~~~~~~~~~\qquad\quad\quad
\leq\tilde{C}\mrm{e}^{-\tilde{\kappa} t}
\Big(1+(\bar{v}_{\frac{\tilde{\kappa} t}{C_{\varepsilon}}}+1)
(x_1+\Delta_1)+\sqrt{x_1+\Delta_1}+|x_2|+\Delta_2\Big),
 \eeqnn
where
 $$
\Delta_1=\int_{\mbb{G}}y_1\,\pi(\mrm{d}y),
\quad \Delta_2=\int_{\mbb{G}}|y_2|\,\pi(\mrm{d}y).
 $$
In order to prove the exponential ergodic property, it suffices to prove $\Delta_1$ and $\Delta_2$
are finite. For $u\in U,$ setting $u=(u_1,0),u_1\leq0,$ then
 $$
\int_{\mbb{G}}e^{u_1y_1}\,\pi(\mrm{d}y)=\exp\Big\{\int_0^{\infty}P(V_1(s,u))\,\mrm{d}s\Big\},
\quad \frac{\partial V_1(s,u)}{\partial s}=\tilde{\phi}_0(V_1(s,u)),
 $$
where $P$ and $\tilde{\phi}_0$ are two functions defined on $\mbb{R}_{-}$ such that
 $$
P(x)=a_2x+\int_{\mbb{G}}(\mrm{e}^{xz_1}-1)\,n(\mrm{d}z),
\quad \tilde{\phi}_0(x)=-a_1x+(\alpha_{11}+\alpha_{12})x^2+\int_{\mbb{G}}
(\mrm{e}^{xz_1}-1-xz_1)\,m(\mrm{d}z).
 $$
 Following the proof of Theorem 10.4 in Li (2020a), in the case of $a_1>0,$ we can see that for $u\in U$,
 $\lim\limits_{s\to \infty}V_1(s,u)=0$ and
  $$
 \int_0^t P(V_1(s,u))\,\mrm{d}s
 =-\int_{V_1(t,u)}^{u_1}\frac{P(z)}{\tilde{\phi}_0(z)}\,\mrm{d}z.
  $$
It implies that
 $$
\int_{\mbb{G}}e^{u_1y_1}\,\pi(\mrm{d}y)=\exp\Big\{-\int_{0}^{u_1}\frac{P(z)}{\tilde{\phi}_0(z)}\,\mrm{d}z\Big\}.
 $$
By differentiating both sides of above at $\lambda_1=0$ we get $\Delta_1=a_1^{-1}[a_2+\int_{\mbb{G}}z_1\,n(\mrm{d}z)]<\infty.$ On the other hand, $\{Z_t\}_{t\geq0}$ can be constructed by
 $$
Z_t=\mrm{e}^{-b_2t}Z_0+\int_0^t\mrm{e}^{-b_2(t-s)}\mrm{d}L^Y_s,\quad t\geq0,
 $$
where
 \beqnn
&&\mrm{d}L^Y_t=-(b_0+b_1Y_t)\mrm{d}t+\sigma\int_0^1\,W_0(\mrm{d}t,\mrm{d}u)
+\sqrt{2\alpha_{21}}\int_0^{Y_t}\,W_1(\mrm{d}t,\mrm{d}u)\\
&&~~~~~~
+\sqrt{2\alpha_{22}}\int_0^{Y_t}\,W_2(\mrm{d}t,\mrm{d}u)
+\int_{\mbb{G}} z_2\,\tilde{N}(\mrm{d}t,\mrm{d}z)
+\int_0^{Y_{t-}}\int_{\mbb{G}}z_2\,\tilde{M}(\mrm{d}t,\mrm{d}u,\mrm{d}z).
 \eeqnn
We assume that $Z_0=y_2,Y_0=y_1.$ By Cauchy-Schwartz inequality and the result for the first
moment of CBI-processes; see, e.g., Li (2020a), pp.33, one can see that
 \beqnn
&&\mbb{E}|\int_0^t \mrm{e}^{b_2s}\,\mrm{d}L^Y_s|
\leq|b_0|\frac{\mrm{e}^{b_2t}-1}{b_2}+|b_1|\int_0^t\mrm{e}^{b_2s}\Big(y_1\mrm{e}^{-a_1s}
+\frac{\gamma}{a_1}(1-\mrm{e}^{-a_1s})\Big)\,\mrm{d}s\\
&&~~~~~~~~~~~~~~~~~~~~~
+\Big[\int_{\mbb{G}}z_2^2\,m(\mrm{d}z)\int_0^t\mrm{e}^{2b_2s}\Big(y_1\mrm{e}^{-a_1s}
+\frac{\gamma}{a_1}(1-\mrm{e}^{-a_1s})\Big)\,\mrm{d}s\Big]^{\frac{1}{2}}\\
&&~~~~~~~~~~~~~~~~~~~~~
+\Big[\int_{\mbb{G}}z_2^2\,n(\mrm{d}z)(\mrm{e}^{b_2t}-1)b^{-1}_2\Big]^{\frac{1}{2}},
 \eeqnn
where $\gamma=a_2+\int_{\mbb{G}}z_1\,n(\mrm{d}z)$. Then
 $$
\sup_{t>0}\mbb{E}|\int_0^t\mrm{e}^{-b_2(t-s)}\,\mrm{d}L^Y_s|
\leq\frac{2|b_1|}{b_2-a_1}y_1+C_5+\Big[\frac{2\int_{\mbb{G}}z_2^2\,m(\mrm{d}z)}
{2b_2-a_1}y_1+C_6\Big]^{\frac{1}{2}}
 $$
for some constants $C_5$ and $C_6$. Further, for any $k\geq1$ and $t>0$,
 $$
\mbb{E}(|Z_t(y)|\wedge k)\leq\mrm{e}^{-b_2t}|y_2|\wedge k
+\sup_{t>0}\mbb{E}|\int_0^t\mrm{e}^{-b_2(t-s)}\,\mrm{d}L^Y_s|,
 $$
and so
 $$
\int_{\mbb{G}}(|y_2|\wedge k)\,\pi(\mrm{d}y)\leq\int_{\mbb{G}}(\mrm{e}^{-b_2t}|y_2|\wedge k)\,
\pi(\mrm{d}y)+C_7, \quad t>0,\,\,k\geq1
 $$
for some constant $C_7$ since $\Delta_1<\infty$. Letting first $t\to \infty$ and then
$k\to \infty$, it follows from dominated convergence theorem and $b_2>0$ that
  $$
 \int_{\mbb{G}}|y_2|\,\pi(\mrm{d}y)<\infty.
  $$
That completes the proof.
\qed

\bigskip

\section{Strong Feller property}

\setcounter{equation}{0}

In this section, we study the strong Feller property of (1+1)-dimensional affine Markov processes constructed by (\ref{3.2})--(\ref{3.3}). We first formulate the following condition:

\noindent$ (\textbf{D})$ There exists a nonnegative measurable function $\varrho_0$ on $\mbb{R}$
such that
 $$
n(\mbb{R}_{+}\times \mrm{d}z_2)
\geq\sigma_0(\mrm{d}z_2)
:=\varrho_0(z_2)\mrm{d}z_2,
\quad \sigma_0(\mbb{R})>0.
 $$
For a positive integrally function $\mrm{g}$ defined on $\mbb{R}$, for $k\geq1,$ let $\sigma_k:=(k\mrm{g}\wedge\rho_0)(z_2)\,\mrm{d}z_2,$ $\sigma_k(\mbb{R}):=\int_{\mbb{R}}(k\mrm{g}\wedge\rho_0)(z_2)\,\mrm{d}z_2$. There exists $K\geq1$ such that for all $k\geq K$, the measure $\sigma_k$ satisfies
 \beqnn
\limsup\limits_{\rho\to 0}\Big[\frac{\sup_{|a|\leq\rho}\|\sigma_k
-\delta_{a}\ast\sigma_k\|_{\mrm{var}}}{\rho}\Big]<\infty,
\quad \int_{\mbb{R}}|z_2|\,
\sigma_k(\mrm{d}z_2)<\infty.
 \eeqnn

Under Condition $\textbf{(D)}$, $\sigma_k(\mbb{R})<\infty$ for $k\geq1$. Then for $k\geq K$ we can define a compound Poisson process as follows:
 $$
L_t^k:=\int_0^t\int_{\mbb{G}}z_2\,N_{k}(\mrm{d}s,\mrm{d}z),
 $$
where $N_{k}(\mrm{d}s,\mrm{d}z)$ is a Poisson random measure with intensity $\mrm{d}sn_k(\mrm{d}z)$ satisfying $n_{k}(\mbb{R}_{+}\times\mrm{d}z_2)=\sigma_k(\mrm{d}z_2)$.
We define a sequence of stopping times $\Upsilon^k_n=\inf\{t>\Upsilon^k_{n-1}:L^k_t\neq L_{t-}^k\}$ with convention $\Upsilon^k_0=0$. For $i\geq1$ let $\tau^k_i=\Upsilon^k_{i}-\Upsilon^k_{i-1}$ and $U^{k}_i= \int_{\{\Upsilon^k_{i}\}} \int_{\mbb{G}} z_2\, N_k(\mrm{d}s,\mrm{d}z).$
Then it is easy to see that $(\tau^k_i)_{i\geq1}$ are i.i.d. random variables which are exponentially distributed with intensity $\sigma_k(\mbb{R})$ and $(U^k_i)_{i\geq1}$ are i.i.d. random variables on $\mbb{R}$ with distribution $\sigma_k/\sigma_k(\mbb{R})$. Moreover, the two sequences are independent of each other. Let $N_t^k=\sup\{j\geq1:\sum_{i=1}^j\tau_i^k\leq t\}$. Then $(N_t^k)_{t\geq0}$ is a Poisson process of intensity $\sigma_k(\mbb{R})$. Now we can rewrite
 $$
\int_0^t\int_{\mbb{G}}z_2\,N_k(\mrm{d}s,\mrm{d}z)=\sum\limits_{i=1}^{N^k_t}U^k_i
 $$
with $\sum_{i=1}^0=0$ by convention.

\bglemma\label{t5.1} Suppose Condition $\textbf{(D)}$ holds. For $k\geq K$ let $P^k_t$ be defined by
 $$
P^k_tf(x):=\mbb{E}\Big[f(Y_t(x_1),Z_t(x))\mbf{1}_{\{\tau^k_1\leq t\}}\Big],
\quad f\in B_b(\mbb{G}),\quad x\in \mbb{G}.
 $$
Then there exist a constant $C_{k,8}>0$ and non-negative functions $t\mapsto C_1(t)$ and $t\mapsto C_2(t)$
such that for any given $x,y\in\mbb{G}$,
 \beqnn
&&\sup_{\|f\|\leq1}|P^k_tf(x)-P^k_tf(y)|\\
&&\leq 2\mrm{e}^{-\sigma_k(\mbb{R})t}+C_{k,8}\Big\{|x_2-y_2|+\sqrt{C_1(t)|x_1-y_1|}
+C_2(t)|x_1-y_1|\Big\}.
 \eeqnn
\edlemma

\proof
By a modified proof in Proposition \ref{t4.2}, we have
 \beqnn
&&\sup_{\|f\|\leq1}|P^k_tf(x)-P^k_tf(y)|\\
&&\leq2\mrm{e}^{-\sigma_k(\mbb{R})t}+\frac{\Lambda_k}{\sigma_k(\mbb{R})}
\Big\{|x_2-y_2|+\sqrt{C_1(t)|x_1-y_1|}+C_2(t)|x_1-y_1|\Big\},
 \eeqnn
where
 \beqnn
&&C_1(t)=\bar{c}\cdot\frac{\mrm{e}^{(2b_2-a_1)t}-1}{2b_2-a_1},\quad
C_2(t)=\frac{|b_1|(\mrm{e}^{(b_2-a_1)t}-1)}{b_2-a_1},\\
&&\bar{c}=4\max\Big\{8(\alpha_{21}\vee\alpha_{22}),\int_{\mbb{G}}|z_2|^2\,m(\mrm{d}z)\Big\},\quad
\Lambda_k=\sup\limits_{\rho>0}\Big[\frac{\sup_{|a|\leq\rho}\|\sigma_k-\delta_{a}\ast \sigma_k\|_{\mrm{var}}}{\rho}\Big]
 \eeqnn
by conventions $C_1(t)=\bar{c}t$ for $2b_2=a_1$ and $C_2(t)=|b_1|t$ for $b_2=a_1.$ And we prove the
required assertion by setting
$C_{k,8}=\frac{\Lambda_k}{\sigma_k(\mbb{R})}.$
\qed

\bgtheorem\label{t5.2}
Suppose Conditions $\textbf{(A,D)}$ hold with $\sigma_0(\mbb{R})=\infty$. Then the affine Markov process $(X_t)_{t\geq0}$ satisfies the strong Feller property.
\edtheorem

\proof
It follows from (\ref{3.4}) that for $x,y\in\mbb{G}$ we have
 \beqnn
\sup_{\|f\|\leq1}|P_tf(x)-P_tf(y)|\leq2\bar{v}_t|x_1-y_1|
+\sup_{\|f\|\leq1}|P^k_tf(x)-P^k_tf(y)|,\quad k\geq K,
 \eeqnn
which together with Lemma \ref{t5.1} imply that
 \beqnn
\lim\limits_{y\to x}\sup_{\|f\|\leq1}|P_tf(x)-P_tf(y)|
\leq2\mrm{e}^{-\sigma_k(\mbb{R})t},
\quad k\geq K,\quad t>0.
 \eeqnn
Since $\sigma_k(\mbb{R})\uparrow\infty$ as $k\uparrow\infty,$ we conclude the strong Feller property of
$(X_t)_{t\geq0}$.
\qed

\end{document}